\documentclass[final]{siamltex}

\usepackage{graphicx}
\usepackage{pstricks}
\usepackage{psfig}

\title{Finite difference schemes as a matrix equation}

\author{Claire David \thanks{Universit\'e Pierre et Marie Curie-Paris 6,
Laboratoire de Mod\'elisation en M\'ecanique, UMR CNRS 7607, Bo\^ite
courrier $n^0162$, 4 place Jussieu, 75252 Paris cedex 05, France
({\tt david@lmm.jussieu.fr}).}}

\begin{document}

\maketitle

\begin{abstract}

Finite difference schemes are here solved by means of a linear
matrix equation. The theoretical study of the related algebraic
system is exposed, and enables us to minimize the error due to a
finite difference approximation.
\end{abstract}

\begin{keywords}
Finite difference schemes, Sylvester equation
\end{keywords}

\begin{AMS}
15A15, 15A09, 15A23
\end{AMS}

\pagestyle{myheadings} \thispagestyle{plain} \markboth{CL.
DAVID}{Finite difference schemes as a matrix equation}

\section{Introduction: Scheme classes}
\label{sec:intro}

\indent Finite difference schemes used to approximate linear
differential equations are usually solved by means of a recursive
calculus. We here propose a completely different approach, which
uses an equivalent matrix equation.\\
\\

\noindent Consider the transport equation:
\begin{equation}
\label{transp} \frac{\partial u}{\partial t}+c\,\frac{\partial
u}{\partial x}=0 \,\,\, , \,\,\, x \,\in \,[0,L], \,\,\,t \,\in
\,[0,T]
\end{equation}

\noindent with the initial condition $u(x,t=0)=u_0(x)$.

\bigskip

\begin{proposition}

\noindent A finite difference scheme for this equation can be
written under the form:
\begin{equation} \label{scheme} {{{{{\alpha \, u}}_i}}^{n+1}}+
   {{{{{\beta \,u}}_i}}^{n}}
    +{{{{{\gamma \,u}}_i}}^{n-1}}
      +\delta \,{{{u_{i+1}}}^n}+{{{{{\varepsilon \, u}}_{i-1}}}^n}
            +{{{{{\zeta \,u}}_{i+1}}}^{n+1}}
              +{{{{{\eta \,u}}_{i-1}}}^{n-1}}+{{{{{\theta \,u}}_{i-1}}}^{n+1}}+\vartheta \,{{ u}_{i+1}}^{n-1} =0
              \end{equation}

\noindent where:
\begin{equation}
{u_l}^m=u\,(l\,h, m\,\tau)
\end{equation}
\noindent  $l\, \in \, \{i-1,\, i, \, i+1\}$, $m \, \in \, \{n-1,\,
n, \, n+1\}$, $j=0, \, ..., \, n_x$, $n=0, \, ..., \, n_t$, $h$,
$\tau$ denoting respectively the mesh size and time step ($L=n_x\,h$, $T=n_t\,\tau$).\\
The Courant-Friedrichs-Lewy number ($cfl$) is defined as $\sigma = c \,\tau / h$ .\\
\\

A numerical scheme is  specified by selecting appropriate values of
the coefficients  $\alpha$, $\beta$, $\gamma$, $\delta$,
$\varepsilon$, $\zeta$, $\eta$, $\theta$ and  $\vartheta$ in
equation (\ref{scheme}). Values corresponding to numerical schemes
retained for the present works are given in Table \ref{SchemeTable}.
\end{proposition}
\bigskip

\begin{table}
\caption{Numerical scheme coefficient.}
\begin{center} \footnotesize
{\begin{tabular}{cccccccccc} \hline
Name & $  \alpha $ & $ \beta$ & $\gamma$ & $\delta $  & $\epsilon$ & $\zeta$ & $\eta$ & $\theta$ & $\vartheta$ \\
\hline
%1st order Upwind & $\frac{1}{ \tau} $ & $  \frac{-1}{ \tau}  + \frac{c}{ h} $ & 0 & 0 & $ \frac{ - c}{h}  $ & 0 & 0 & 0 & 0 \\
Leapfrog & $\frac{1}{2 \tau} $ & 0 &  $\frac{-1}{2 \tau} $ & $  \frac{c}{ 2 h}  $ & $   \frac{-c}{ 2 h}  $ & 0 & 0 & 0  & 0 \\
Lax & $\frac{1}{ \tau} $ &  0 & 0 & $  \frac{-1}{ 2 \tau} + \frac{c}{ 2 h } $ & $  \frac{-1}{2 \tau}  - \frac{c}{ 2 h} $ & 0 & 0 & 0  & 0 \\
Lax-Wendroff & $\frac{1}{ \tau} $ & $  \frac{-1}{ \tau}  + \frac{c^2 \tau}{ h ^2} $ & 0 & $ \frac{( 1- \sigma ) c}{ 2 h} $&  $  \frac{-( 1+ \sigma ) c}{ 2 h } $ & 0 & 0 & 0 & 0 \\
Crank-Nicolson & $ \frac{1}{ \tau}  + \frac{c}{ h ^2} $ &  $
\frac{-1}{ \tau}  + \frac{c}{ h ^2} $ & 0 & $ \frac{-c}{ h ^2} $ & $
\frac{-c}{ h ^2 } $ & 0 & $ \frac{-c}{ h ^2} $ & $ \frac{-c}{ h ^2}
$ & 0
\end{tabular}}
\end{center}
\label{SchemeTable}
\end{table}

The number of time steps will be denoted $n_t$, the number of space
steps, $n_x$. In general, $n_x\gg n_t$.\\
\\
%We have choosen here to restrain our study to the case $n_x=n_t$.\\
%We will note: $N=n_x=n_t$.\\

%The case $n_x <
%n_t$ can be easily deduced from it, as it will be exposed later.\\
\bigskip

The paper is organized as follows. The equivalent matrix equation is
exposed in section \ref{Sylv}. Specific properties of the involved
matrices are exposed in section \ref{Properties}. A way of
minimizing the error due to the finite difference approximation is
presented in section \ref{MinErr}.

\section{The Sylvester equation}
\label{Sylv}
\subsection{Matricial form of the finite differences problem}

\noindent Let us introduce the rectangular matrix defined by:
\begin{equation}
U=[{{{u_i}}^n}{]_{\, 1\leq i\leq {n_x-1}, \, 1\leq n\leq {n_t} \, }}
\end{equation}

\begin{theorem}
The problem (\ref{scheme}) can be written under the following
matricial form:
\begin{equation}
\label{Eq} {M_1}\,U +U\,M_2+{\cal{L}}(U)=M_0
\end{equation}

\noindent where $M_1$ and $M_2$ are square matrices respectively
$n_x-1$ by $n_x-1$, $n_t$ by $n_t$, given by:
\begin{equation}
\begin{array}{ccc}% autant de c que de colonnes
{M_1}= \left (
\begin{array}{ccccc}% autant de c que de colonnes
\beta & \delta & 0 & \ldots &  0 \\
\varepsilon& \beta & \ddots & \ddots &  \vdots \\
0 &  \ddots & \ddots & \ddots &  0\\
\vdots &  \ddots & \ddots & \beta &  \delta\\
0 &  \ldots & 0 & \varepsilon &  \beta\\
\end{array} \right ) &
 & {M_2}= \left (
\begin{array}{ccccc}% autant de c que de colonnes
0 & \gamma& 0 & \ldots &  0 \\
\alpha & 0 & \ddots & \ddots &  \vdots \\
0 &  \ddots & \ddots & \ddots &  0\\
\vdots &  \ddots & \ddots & \ddots &  \gamma\\
0 &  \ldots & 0 & \alpha &  0\\
\end{array} \right )
\end{array}
  \end{equation}

\noindent the matrix $M_0$ being given by:
\bigskip
  \begin{equation}
  \label{M0}
\scriptsize{{M_0}= \left (
\begin{array}{ccccc}% autant de c que de colonnes
    -\gamma \,u_1^0
      -\varepsilon \,u_{0}^1-\eta\,u_0^0-\theta\,u_{0}^{2}-\vartheta\,u_{2}^{0}
            &   -\varepsilon \,u_{0}^2 -\eta\,u_0^1-\theta\,u_{0}^{3}& \ldots & \ldots   &
            -\varepsilon \,u_{0}^{n_t}-\eta\,u_{0}^{n_t-1} \\
 -\gamma \,u_2^0-\eta\,u_{1}^{0}-\vartheta\,u_{3}^{0}
       &  0
 & \ldots  & \ldots    &  0 \\
\vdots  &  \vdots & \vdots & \vdots &  \vdots\\
 -\gamma \,u_{n_x-2}^0-\eta\,u_{n_x-2}^{0}-\vartheta\,u_{n_x-1}^{0}
       &  0
 & \ldots  & \ldots   &  0 \\
    -\gamma \,u_{n_x-1}^0-\delta \,u_{n_x}^{1}-\eta\,u_{n_x-2}^{0}-\zeta\,u_{n_x}^{2}-\vartheta\,u_{n_x}^{0}
      &   -\delta \,u_{n_x}^{2} -\zeta\,u_{n_x}^{3}-\vartheta\,u_{n_x}^{1}& \ldots
       & \ldots &  -\delta \,u_{n_x}^{n_t}-\vartheta\,u_{n_x}^{n_t-1}\\
\end{array} \right )}
  \end{equation}

\bigskip
\noindent and where ${\cal {L}}$  is a linear matricial operator
which can be written as: \begin{equation} {\cal {L}}={\cal
{L}}_1+{\cal {L}}_2+{\cal {L}}_3+{\cal {L}}_4
\end{equation}
\noindent where ${\cal {L}}_1$, ${\cal {L}}_2$, ${\cal {L}}_3$ and
${\cal {L}}_4$ are given by:
\begin{equation}
\begin{array}{ccc}% autant de c que de colonnes
 {\cal {L}}_1(U)  = \zeta  \left (
\begin{array}{ccccc}% autant de c que de colonnes
u_2^2 & u_2^3 & \ldots & u_2^{n_t} &  0\\
u_3^2 & u_3^3 & \ldots & \vdots &   \vdots \\
\vdots &  \vdots & \ddots & \vdots &  \vdots\\
u_{n_x-1}^2 & u_{n_x-1}^3  & \ldots &  u_{n_x-1}^{n_t}& 0  \\
0 & 0  & \ldots &  0& 0  \\
\end{array}
 \right )
  &
  &  {\cal {L}}_2(U) = \eta \left (
\begin{array}{ccccc}% autant de c que de colonnes
0 & 0  & \ldots &  0& 0 \\
0 & u_1^1  & u_1^2 & \ldots &  u_1^{n_t-1} \\
0 & u_1^0  & u_1^1 & \ldots &  u_2^{n_t-1} \\
\vdots &  \vdots & \vdots &  \ddots & \vdots \\
0 & u_{n_x-2}^1  & u_{n_x-2}^2 & \ldots &  u_{n_x-2}^{n_t-1} \\
\end{array}
 \right )
\end{array}
  \end{equation}

\begin{equation}
\begin{array}{ccc}% autant de c que de colonnes
{\cal {L}}_3(U) = \theta \left (
\begin{array}{ccccc}% autant de c que de colonnes
0 & \ldots &  \ldots &    \ldots &  0 \\
u_1^2 & u_1^3  & \ldots &  u_1^{n_t}& 0 \\
u_2^2 & u_2^3  & \ldots &  u_2^{n_t}& 0 \\
\vdots &  \vdots & \vdots &  \vdots & \vdots \\
u_{n_x-2}^2 & u_{n_x-2}^3  & \ldots &  u_{n_x-2}^{n_t}& 0 \\
\end{array}
 \right )
  &
  &  {\cal {L}}_4(U) = \vartheta \left (
\begin{array}{ccccc}% autant de c que de colonnes
0 & u_2^1 & u_2^2 & \ldots  &  u_2^{n_t-1} \\
0 & u_3^1 & u_3^2 & \ldots  &  u_3^{n_t-1} \\
\vdots &  \vdots & \ddots  & \ddots &  \vdots\\
0 & u_{n_x-1}^1& \ldots & \ldots  &  u_{n_x-1}^{n_t-1} \\
0 &  0 & \ldots &   \ldots &  0\\
\end{array}
 \right )
 \end{array}
 \end{equation}
\bigskip

\end{theorem}

\bigskip

\begin{proposition}
\noindent The second member matrix $M_0$ bears the initial
conditions, given for the specific value $n=0$, which correspond to
the initialization process when computing loops, and the boundary
conditions, given for the specific values $i=0$, $i=n_x$.
\end{proposition}

\bigskip

\noindent Denote by $u_{exact}$ the exact solution of (\ref{transp}).\\
\noindent The corresponding matrix $U_{exact}$ will be:

\begin{equation}
U_{exact}=[{{{U_{{exact}_i}}}^n}{]_{\, 1\leq i\leq {n_x-1},\, 1\leq
n\leq {n_t}\, }} \end{equation} where:

\begin{equation}
{U_{exact}}_i^n=U_{exact}(x_i,t_n)
 \end{equation}

\noindent with $x_i=i \; h$, $t_n=n \; \tau$. \bigskip

\bigskip

  \begin{definition}
\noindent We will call \textit{error matrix} the matrix defined by:
 \begin{equation}
 \label{err}
E=U-U_{exact}
   \end{equation}
  \end{definition}

\bigskip

\noindent Consider the matrix $F$ defined by:
\begin{equation} F={M_1}\,U_{exact}+U_{exact}\,M_2 + {\cal{L}}(U_{exact})-M_0\end{equation}

\bigskip

  \begin{proposition}
\noindent The \textit{error matrix} $E$ satisfies:

   \begin{equation}
   \label{eqmtr}
{M_1}\,E+E\,M_2+{\cal{L}}(E)=F
   \end{equation}

 \end{proposition}

\newpage

\subsection{The matrix equation}

\subsubsection{Theoretical formulation}

   \begin{theorem}
\noindent Minimizing the error due to the approximation induced by
the numerical scheme is equivalent to minimizing the norm of the
matrices $E$ satisfying (\ref{eqmtr}).
   \end{theorem}

   \bigskip

{\em Note:} \noindent Since the linear matricial operator
${\cal{L}}$ appears only in the Crank-Nicolson scheme, we will
restrain our study to the case ${\cal{L}}=0$. The generalization to
the case ${\cal{L}} \neq 0$ can be easily deduced.

   \bigskip

   \begin{proposition}

   \noindent The problem is then the determination of the minimum norm solution
of:

   \begin{equation}
   \label{SylvErr}
{M_1}\,E+E\,M_2=F
   \end{equation}

\noindent which is a specific form of the Sylvester equation:

   \begin{equation}
   \label{SylvGen}
AX+XB=C
   \end{equation}
where $A$ and $B$ are respectively $m$ by $m$ and $n$ by $n$
matrices, $C$ and $X$, $m$ by $n$ matrices.

   \end{proposition}

\bigskip

 \noindent The solving of the Sylvester
equation is generally based on Schur decomposition: for a given
square $n$ by $n$ matrix $A$, $n$ being an even number of the form
$n=2\,p$, there exists a unitary matrix $U$ and a upper triangular
block matrix $T$ such that:
\begin{equation}
A = U^*TU
   \end{equation}
\noindent where $U^*$ denotes the (complex) conjugate matrix of the
transposed matrix $^T U$. The diagonal blocks of the matrix $T$
correspond to the complex eigenvalues $\lambda_i$ of $A$:
\begin{equation}
T = \left (
\begin{array}{ccccc}% autant de c que de colonnes
T_1 & 0 & \ldots  & \ldots   & 0 \\
0 & \ddots & \ddots   & \ddots &  \vdots \\
\vdots &  \ddots & T_i  & \ddots &  \vdots \\
\vdots &  \ddots & \ddots &   \ddots  & 0 \\
0 &  0 & \ldots  & 0 & T_p\\
\end{array}
 \right )
   \end{equation}
\noindent where the block matrices $T_i$, $i=1,\ ...,\, p$ are given
by:
\begin{equation}
\left (\begin{array}{cc}% autant de c que de colonnes
\mathcal{R}e\, [\lambda_i] & \mathcal{I}m\, [\lambda_i]  \\
-\,\mathcal{I}m\, [\lambda_i] & \mathcal{R}e\, [\lambda_i]  \\
\end{array}
 \right )
   \end{equation}

\noindent $\mathcal{R}e$ being the real part of a complex number,
and $\mathcal{I}m$ the imaginary one.\\
\bigskip \noindent Due to this decomposition, the Sylvester equation
require, to be solved, that the dimensions of the matrices be even
numbers. We will therefore, in the following, restrain our study to
$n_x-1$ and $n_t$ being even numbers. So far, it is interesting to
note that the Schur decomposition being more stable for higher order
matrices, it perfectly fits finite difference problems.

\bigskip

\noindent Complete parametric solutions of the generalized Sylvester
equation (\ref{SylvErr}) is given in \cite{Berman}, \cite{Gail}.\\

\noindent As for the determination of the solution of the Sylvester
equation, it is a major topic in control theory, and has been the
subject of numerous works (see \cite{Van}, \cite{Hearon},
\cite{Tsui}, \cite{Zhou},
\cite{Duan}, \cite{Duan2}, \cite{Kir}).\\
\noindent In \cite{Van}, the method is based on the reduction of the
he observable pair $(A,C)$ to an observer-Hessenberg  pair $(H, D)$,
$H$ being a block upper Hessenberg matrix. The reduction to the
observer-Hessenberg form $(H, D)$ is achieved by means of the
staircase algorithm
(see \cite{Bol}, ...). \\
\noindent In \cite{Zhou}, in the specific case of $B$ being a
companion form matrix, the authors propose a very neat general
complete parametric solution, which is expressed in terms of the
controllability of the matrix pair $(A,B)$, a symmetric matrix
operator, and a parametric matrix in the Hankel form.\\
We recall that a companion form, or Frobenius matrix is one of the
following kind:
 \begin{equation}
B = \left (
\begin{array}{cccccc}% autant de c que de colonnes
0 & \ldots & \ldots  & \ldots  & 0 &  -b_{0} \\
1 & 0 & \ldots  & \ldots & 0 &  -b_{1} \\
0 &  1 & 0  & \ldots & \vdots &  \vdots \\
\vdots &  0 & \ddots & \vdots &  \vdots & \vdots \\
0 &  0 & \ldots & 1 & 0 &  -b_{p-1}\\
\end{array}
 \right )
   \end{equation}

\noindent These results can be generalized through matrix block
decomposition to a block companion form matrix:
\begin{equation}
M_2 =  \left (
\begin{array}{ccccc}% autant de c que de colonnes
{M_2^{B}}^1 & 0 & \ldots  & \ldots   &  0 \\
0 &  {M_2^{B}}^2 & 0 &  \ldots &  0 \\
0 &  0 & \ddots  & \ddots  &  \vdots \\
\vdots &  0 & \ddots & \ddots  & 0 \\
0 &  0 &  \ldots & 0 &  {M_2^{B}}^k\\
\end{array}
 \right )
   \end{equation}

 \noindent the ${M_2^{B}}^p$, $1\leq p \leq k$ being companion
 form matrices.\\
\\
\\
\noindent Another method is presented in \cite{Var}, where the
determination of the minimum-norm solution
of a Sylvester equation is specifically developed.\\
\noindent The accuracy and computational stability of the solutions
is examined in \cite{Deif}.
\bigskip

\bigskip

\subsubsection{Existence condition of the solution }

\begin{theorem} In the specific cases of the Lax and Lax-Wendroff schemes, (\ref{SylvGen}) has a unique
solution.

\end{theorem}

\begin{proof} \noindent (\ref{SylvGen}) has a unique
solution if and only if $A$ and $B$ have no common eigenvalues.\\
\noindent The characteristic polynomials $P_{M_1}$, $P_{M_2}$ of
$M_1$ and $M_2$, can be respectively calculated as the determinants
of respectively $\frac{n_x-1}{2}$
 by $\frac{n_x-1}{2}$,  $\frac{n_t}{2}$
 by $\frac{n_t}{2}$ diagonal block matrices:

\begin{equation}
P_{M_1}(\lambda)=((\lambda-\beta)^2-\delta\,\varepsilon)^{\frac{n_x-1}{2}}\,\,\,,
\,\,\, P_{M_2}(\lambda)=(\lambda^2-\alpha\,\gamma)^{\frac{n_t}{2}}
   \end{equation}
\noindent In the specific cases of the Lax and Lax-Wendroff schemes:
$\alpha\,\gamma=0\neq \pm \beta +\sqrt{\delta \varepsilon}$.
    Hence, (\ref{eqmtr}) has a
unique solution, which accounts for the consistency of the given
problem.

\end{proof}

\begin{corollary} In the specific case of the Leapfrog scheme, (\ref{SylvGen}) has a unique
solution if and only if $\tau \neq \frac{h}{c}$.

\end{corollary}

\section{Specific properties of the matrices $M_1$ and $M_2$}

\label{Properties}

\subsubsection{Inversibility of the matrix $M_1$ }

\begin{theorem} In the specific case of the Lax and Leapfrog schemes, $M_1$ is inversible.

\end{theorem}

\begin{theorem} In the specific case of the Lax-Wendroff scheme, $M_1$ is inversible if and only if

\begin{equation}
\big(\frac{-1}{ \tau}  + \frac{c^2 \tau}{ h ^2}\big)^2  \neq \frac{(
 \sigma^2-1 ) c^2}{ 4\, h^2}
\end{equation}
\end{theorem}

\begin{proof} \noindent The determinant $D_1$ of $M_1$ can be calculated as the determinant
of a $\frac{n_x-1}{2}$
 by $\frac{n_x-1}{2}$ diagonal block matrix:

\begin{equation}
D_1=(\beta^2-\delta\,\varepsilon)^{\frac{n_x-1}{2}}
\end{equation}

\end{proof}

\subsubsection{Nilpotent components of the matrix $M_2$}

 \begin{proposition}
\noindent $M_2$ can be written as:
 \begin{equation}
M_2=\alpha\ ^T N \,+\,\gamma\,  N
   \end{equation}
\noindent where $N$ is the nilpotent matrix; $^T N $ denotes the
corresponding transposed matrix:
 \begin{equation}
N= \left (
\begin{array}{ccccc}% autant de c que de colonnes
0 & 1 & 0  & \ldots &  0 \\
0 & 0 & \ddots  & \ddots  &  \vdots \\
\vdots  &  \vdots & \ddots & \ddots &  0\\
\vdots &  \vdots & \ldots & 0 & 1\\
0 &  0 & \ldots & \ldots & 0\\
\end{array} \right )
  \end{equation}
 \end{proposition}

\bigskip

 \begin{proposition}

 \noindent
\noindent For either $\alpha=0$, or $\gamma=0$, $M_2$ will thus be
nilpotent, of order $n_t$.

 \end{proposition}

\bigskip

 \begin{proposition}

 \noindent For $\gamma=0$, if $M_1$ is inversible, the solution at $t=n_t\,dt$ can then be immediately
 determined.

 \end{proposition}

\bigskip

 \begin{proof}

 \noindent In such a case, multiplying (\ref{Eq}) on the right side by ${M_2}^{n_t-1}$
leads to:
 \begin{equation}
M_1 \,U\,{M_2}^{n_t-1}+U\,{M_2}^{n_t}=M_0\,{M_2}^{n_t-1}
   \end{equation}
\noindent i. e.:
\begin{equation}
M_1 \,U\,{M_2}^{n_t-1}=M_0\,{M_2}^{n_t-1}
   \end{equation}
\noindent which leads to:
\begin{equation}
 U\,{M_2}^{n_t-1}=M_1^{-1}\,M_0\,{M_2}^{n_t-1}
   \end{equation}
\bigskip

\noindent Due to:
\begin{equation}
{M_2}^{n_t-1}= \left (
\begin{array}{ccccc}% autant de c que de colonnes
0&0 & 0 & \ldots &  0 \\
\vdots& 0& \ddots & \ddots &  \vdots \\
\vdots &  \vdots & \ddots & \ddots &  0\\
0 &  \vdots & \ddots & 0 &  0\\
1 &  0 & \ldots  & 0 &  0\\
\end{array} \right )
  \end{equation}

 \noindent we have:
\begin{equation}
U\,{M_2}^{n_t-1}= \left (
\begin{array}{ccccc}% autant de c que de colonnes
u_{1n_t}&0 & 0 & \ldots &  0 \\
u_{2n_t}& 0& \ddots & \ddots &  \vdots \\
\vdots &  \vdots & \ddots & \ddots &  0\\
\vdots &  \vdots  & \ddots & 0 &  0\\
u_{n_xn_t} &  0 & \ldots  & 0 &  0\\
\end{array} \right )
  \end{equation}

 \noindent The solution at $t=n_t\,dt$ can then be immediately
 determined.

\end{proof}

\section{Minimization of the error}
\label{MinErr}

\subsection{Theory}

\noindent Calculation yields:

\footnotesize \noindent \begin{equation}\left \lbrace
\begin{array}{ccc}% autant de c que de colonnes
{M_1}\,^T M_1&=& diag \big
 (\left ( \begin{array}{cc}% autant de c que de colonnes
\beta^2+ \delta^2 & \beta\,(\delta+\varepsilon) \\
\beta\,(\delta+\varepsilon) &\varepsilon^2+ \beta^2 \\
\end{array} \right ),\ldots, \left ( \begin{array}{cc}% autant de c que de colonnes
\beta^2+ \delta^2 & \beta\,(\delta+\varepsilon) \\
\beta\,(\delta+\varepsilon) &\varepsilon^2+ \beta^2 \\
\end{array} \right )
 \big )\\ {M_2}\,^T M_2&=&
 diag \big
 (\left ( \begin{array}{cc}% autant de c que de colonnes
\gamma^2 & 0 \\
0 &\alpha^2\\
\end{array} \right ),\ldots, \left ( \begin{array}{cc}% autant de c que de colonnes
\gamma^2 & 0 \\
0 &\alpha^2 \\
\end{array} \right )
\end{array} \right.
  \end{equation}

\normalsize \noindent The singular values of $M_1$ are the singular
values of the block matrix $\big
 (\left ( \begin{array}{cc}% autant de c que de colonnes
\beta^2+ \delta^2 & \beta\,(\delta+\varepsilon) \\
\beta\,(\delta+\varepsilon) &\varepsilon^2+ \beta^2 \\
\end{array} \right )$, i. e. \begin{equation} \frac{1}{2} \,(2 \beta ^2+\delta
^2+\varepsilon
   ^2-(\delta +\varepsilon ) \,\sqrt{4 \beta ^2+\delta
   ^2+\varepsilon ^2-2 \delta \, \varepsilon })  \end{equation} \noindent of order $\frac {n_x-1}{2}$, and \begin{equation} \frac{1}{2} \,(2 \beta
^2+\delta ^2+\varepsilon
   ^2+(\delta +\varepsilon ) \,\sqrt{4 \beta ^2+\delta
   ^2+\varepsilon ^2-2 \delta \, \varepsilon })  \end{equation} \noindent of order $\frac {n_x-1}{2}$.\\

\noindent The singular values of $M_2$ are $\alpha^2 $, of order
$\frac {n_t}{2}$, and $\gamma^2 $, of order $\frac {n_t}{2}$.

\noindent Consider the singular value decomposition of the matrices
$M_1$ and $M_2$:

\begin{equation}
U_1^T\,M_1\,V_1=\left (
\begin{array}{cc}% autant de c que de colonnes
\widetilde{M_1} &0 \\
0& 0\\
\end{array} \right )
\,\,\, ,  \,\,\, U_2^T\,M_1\,V_2=\left (
\begin{array}{cc}% autant de c que de colonnes
\widetilde{M_2} &0 \\
0& 0\\
\end{array} \right )
  \end{equation}

\noindent where $U_1$, $V_1$, $U_2$, $V_2$, are orthogonal matrices.
\noindent $\widetilde{M_1}$, $\widetilde{M_2}$ are diagonal
matrices, the diagonal terms of which are respectively the nonzero
eigenvalues of the symmetric matrices $M_1\,^T M_1$, $M_2\,^T
M_2$.\\

\noindent Multiplying respectively \ref{SylvErr} on the left side by
$^T U_1$, on the right side by $V_2$, yields:

\begin{equation}
U_1^T\,M_1\,E\,V_2+U_1^T\,E\,M_2\,V_2=U_1^T\,F\,V_2
  \end{equation}

\noindent which can also be taken as:
\begin{equation}
^T U_1\,M_1\,V_1 \,^T V_1\,E\,V_2+^T U_1\,E\,^T U_2\,^T
U_2\,M_2\,V_2=U_1^T\,F\,V_2
  \end{equation}

\noindent Set:

\begin{equation}
^T V_1\,E\,V_2= \left (
\begin{array}{cc}% autant de c que de colonnes
\widetilde{E_{11}} & \widetilde{E_{12}} \\
\widetilde{E_{21}} & \widetilde{E_{22}}\\
\end{array} \right )
\,,\,^T U_1\,E\,^T U_2= \left (\begin{array}{cc}% autant de c que de colonnes
\widetilde{\widetilde{E_{11}}} & \widetilde{\widetilde{E_{12}}} \\
\widetilde{\widetilde{E_{21}}} & \widetilde{\widetilde{E_{22}}}\\
\end{array} \right )
  \end{equation}

\begin{equation}
\label{Ftilde} ^T U_1\,F\,V_2= \left (
\begin{array}{cc}% autant de c que de colonnes
\widetilde{F_{11}} & \widetilde{F_{12}} \\
\widetilde{F_{21}} & \widetilde{F_{22}}\\
\end{array} \right )
  \end{equation}
\noindent We have thus:
\begin{equation}\left (\begin{array}{cc}% autant de c que de colonnes
\widetilde{M_{1}}\,\widetilde{E_{11}} & \widetilde{M_{1}}\,\widetilde{E_{12}} \\
0 & 0\\
\end{array} \right )+\left (\begin{array}{cc}% autant de c que de colonnes
\widetilde{\widetilde{E_{11}}}\,\widetilde{M_{2}} & 0 \\
\widetilde{\widetilde{E_{21}}}\,\widetilde{M_{2}} & 0\\
\end{array} \right )=\left (\begin{array}{cc}% autant de c que de colonnes
\widetilde{F_{11}} & \widetilde{F_{12}} \\
\widetilde{F_{21}} & \widetilde{F_{22}}\\
\end{array} \right )
  \end{equation}

\noindent It yields:
\begin{equation}
\left \lbrace
\begin{array}{ccc}% autant de c que de colonnes
\widetilde{M_{1}}\,\widetilde{E_{11}} +
\widetilde{\widetilde{E_{11}}}\,\widetilde{M_{2}}&=&\widetilde{F_{11}}\\
\widetilde{M_{1}} \,\widetilde{E_{12}}&=&\widetilde{F_{12}}\\
\widetilde{\widetilde{E_{21}}}\,\widetilde{M_{2}}&=&\widetilde{F_{21}}\\
\end{array} \right.
  \end{equation}

\noindent One easily deduces:
\begin{equation}
\left \lbrace
\begin{array}{ccc}% autant de c que de colonnes
\widetilde{E_{12}}&=&{\widetilde{M}_{1}}^{-1}\,\widetilde{F_{12}}\\
{\widetilde{\widetilde{E}_{21}}}&=&\widetilde{F_{21}}\,{\widetilde{M_{2}}}^{-1}\\
\end{array} \right.
  \end{equation}

\noindent The problem is then the determination of the
$\widetilde{E_{11}}$ and $\widetilde{\widetilde{E_{11}}}$
satisfying:

\begin{equation}
\label{Pb} \widetilde{M_{1}}\,\widetilde{E_{11}} +
\widetilde{\widetilde{E_{11}}}\,\widetilde{M_{2}}=\widetilde{F_{11}}
  \end{equation}

\noindent Denote respectively by $\widetilde{e_{ij}}$,
$\widetilde{\widetilde{e_{ij}}}$ the components of the matrices
$\widetilde{E}$, $\widetilde{\widetilde{E}}$.\\
\noindent The problem \ref{Pb} uncouples into the independent
problems:\\ \noindent minimize
\begin{equation}
\sum_{i,j} {\widetilde{e_{ij}}}^2+{\widetilde{\widetilde{e_{ij}}}}^2
\end{equation}

 \noindent under the constraint \begin{equation}
\widetilde{M_{1}}_{ii}\,
{\widetilde{e_{ij}}}+\widetilde{M_{2_{ii}}}\,{\widetilde{\widetilde{e_{ij}}}}=\widetilde{F_{11}}_{ij}
  \end{equation}
\noindent This latter problem has the solution:

\begin{equation}\left \lbrace
\begin{array}{ccc}% autant de c que de colonnes
\widetilde{e_{ij}}
&=&\frac{\widetilde{{M_{1}}_{ii}}\,\widetilde{{F_{11}}_{ij}}}
{{\widetilde{{M_{1}}_{ii}}}^2+{\widetilde{{M_{2}}_{jj}}^2}}\\
\widetilde{\widetilde{e_{ij}}} &=&\frac{\widetilde{{M_{2}}
_{jj}}\,\widetilde{{F_{11}}_{ij}}}{{\widetilde{{M_{1}}_{ii}}}^2+{\widetilde{{M_{2}}_{jj}}^2}}\\
\end{array} \right.
  \end{equation}

\noindent The minimum norm solution of \ref{SylvErr} will then be
obtained when the norm of the matrix $\widetilde{{F_{11}}}$ is
minimum.\\
\noindent In the following, the euclidean norm will be considered.

\noindent Due to (\ref{Ftilde}):
 \begin{equation}
  \|\widetilde{{F_{11}}}  \| \leq  \|\widetilde{{F}}  \| \leq \|U_1
  \| \,\|F \|\,\|V_2  \| \leq \|U_1
  \| \,\|V_2  \| \, \|M_1\,U_{exact}+U_{exact}\,M_2-M_0
  \|
   \end{equation}

\noindent $U_1$ and $V_2$ being orthogonal matrices, respectively
$n_x-1$ by $n_x-1$, $n_t$ by $n_t$, we have:

 \begin{equation}
  \|U_1
  \|^2 =n_x-1\,\,\,, \,\,\,\|V_2
  \|^2 =n_t
   \end{equation}

\noindent Also:

 \begin{equation}
  \|M_1
  \|^2 =\frac{n_x-1}{2}\,\big ( 2\,\beta^2+\delta^2+\varepsilon^2 \big )\,\,\,, \,\,\,
  \|M_2
  \|^2 =\frac{n_t}{2}\,\big ( \alpha^2+\gamma^2 \big )
   \end{equation}

\noindent The norm of $M_0$ is obtained thanks to relation
(\ref{M0}).

\noindent This results in:
 \begin{equation}
 \label{Min}
  \|\widetilde{{F_{11}}} \| \leq \sqrt {n_t\,(n_x-1)} \, \left \lbrace  \| U_{exact} \| \,\big (
  \sqrt{\frac{n_x-1}{2}}\,\sqrt{ 2\,\beta^2+\delta^2+\varepsilon^2 }+
   \sqrt{\frac{n_t}{2}}\,\sqrt{\alpha^2+\gamma^2 } \,\big )+
  \|M_0  \|\right \rbrace
   \end{equation}

\noindent  $ \|\widetilde{{F_{11}}} \| $ can be minimized through
the minimization of the right-side member of (\ref{Min}), which is
function of the scheme parameters.

\subsection{Numerical example: the specific case of the Lax scheme}

\noindent For the Lax scheme: $\gamma=0$. The remaining coefficients
(see Table \ref{SchemeTable}) can be normalized through the
following change of variables:
   \begin{equation}
   \left\lbrace
   \begin{array}{rcl}
\overline{\alpha}&=&h\,\alpha \\
\overline{\beta}&=&h\,\beta \\
\overline{\delta}&=&h\,\delta \\
\overline{\varepsilon}&=&h\,\varepsilon \\
\end{array} \right.
   \end{equation}
\noindent Set:
  \begin{equation}
   \overline{M_0}=h\,M_0
   \end{equation}

\noindent Thus:
 \begin{equation}
\|  \overline{M_0}\|^2=\overline{\delta}^2\,\sum_{n=1}^{n_t}
{u_{n_x}^n}^2+\overline{\varepsilon}^2\,\sum_{n=1}^{n_t} {u_0^i}^2
   \end{equation}

\noindent Advect a sinusoidal
signal\\
 \begin{equation} \label{signal} u=\cos\,[\,\frac {2\, \pi}{\lambda} \,(x-c\,t)\,]
\end{equation}
\noindent through the Lax scheme, with Dirichlet boundary
conditions.\\
\noindent The right-side member of (\ref{Min}) is then:
 \begin{equation}
 \left(\frac{1}{2}+\frac{1}{2
\,{cfl}}\right)^2\,{n_t}^2 \,{u_0}^2+\left(\frac{1}{2}-\frac{1}{2
\,{cfl}}\right)^2 \,{n_t}^2
   \,{u_L}^2+\sqrt{\frac{n_x-1}{2}}\,
   \sqrt{\left(\frac{1}{2}+\frac{1}{2
\,{cfl}}\right)^2+\left(\frac{1}{2}-\frac{1}{2
\,{cfl}}\right)^2}+\frac {\sqrt{{n_t}}}{\sqrt{2}\,cfl}
 \end{equation}
 \noindent where $u_0$, $u_L$ respectively denote the Dirichlet
 boundary values at $x=0$, $x=L$.\\
\noindent It is minimal for the admissible value $cfl=1$.

\bigskip

\begin{figure}[ht]
\hspace{3.5cm}
\begin{tabular}{c}
\psfig{height=5cm,width=11cm,angle=0,file=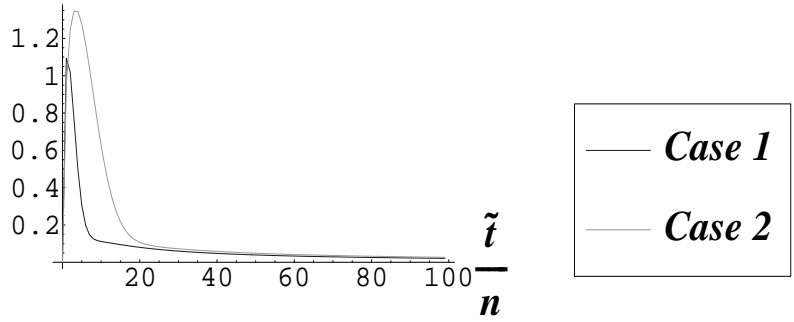}\\
\end{tabular}
\caption{\small{Value of the $L_2$ norm  of the error for different
values of the $cfl$ number}. } \label{Opt}
\end{figure}

\noindent The value of the $L_2$ norm of the error, for two
significant values of the $cfl$ number (Case 1: $cfl=0.9$; Case 2:
$cfl=0.7$), is displayed in Figure \ref{Opt}. The error curve
corresponding to the value of the $cfl$ number closest to 1 is the
minimal one.

\section{Conclusion}

Thanks to the above results, we here propose:
\begin{enumerate}
 \item to study the intrinsic properties of various schemes through those
 of the related matrices $M_1$ and $M_2$;

 \item to optimize finite difference problems through minimization
 of the symbolic expression of the error as a  function of the
 scheme parameters.

 \end{enumerate}

\newpage

\addcontentsline{toc}{section}{\numberline{}References}

\end{document}